\documentclass[12pt]{amsart}

\begin{document}

\begin{center}
{\bf  $L${-}approximation of  $B${-}splines by trigonometric
polynomials }
\end{center}

\vskip .3cm
\begin{center}
A.G. Babenko,  Yu.V. Kryakin
\end{center}
\vskip .3cm

\centerline{\em Dedicated to Professor Eleonora Storozenko on the occasion of
her  eighteen${}^{\text{th}}$ birthday}

\vskip 1cm

\vskip .2cm

Denote by $T_{2n-1}$ a space of real trigonometric  polynomials

 $$
 \tau(x)=\sum_{j=-n+1}^{n-1} \alpha_j\, \exp(2{\pi}i\,jx), \quad
 x\in{\mathbb{T}=\mathbb{R}/\mathbb{Z}},\quad\alpha_{-j}=\overline\alpha_j,
 $$
 and denote by  $T_{2n-1}^\perp$   a
{subspace} of real functions $f\in{L_\infty=L_\infty(\mathbb{T})}$
that are orthogonal to $T_{2n-1}$ with respect to the {``}scalar
product {''}
 $$
 (f,g)=\int_{-1/2}^{1/2}f(u)g(u)\,d{u}=\int_{{\mathbb{T}}}f(u)g(u)\,d{u}.
 $$
\noindent
 It is known that the best approximation in {$L=L(\mathbb{T})$} by
trigonometric polynomials from  $T_{2n-1}$  may be calculate as
follows (see  \cite{n})

 $$
 E_{n}(f)_1:=\inf_{\tau\in T_{2n-1}}\|f-\tau\|_1=
 \sup_{g\in T_{2n-1}^\perp,\ \|g\|=1}\int_{{\mathbb{T}}}g(u)f(u)\,du.
 $$
{Here  $\|\cdot\|_1$ and $\|\cdot\|$  denotes norms in $L$ and
$L_\infty$ respectively.}

This note is a continuation of our papers \cite{bk1, bk2}, devoted
to $L${-}approximation of {$L$-normed} characteristic function $$
{\chi_h(x):=\begin{cases} h^{-1},\quad & x\in{(}-h/2,h/2{)},\\
 0, \quad & x \notin{(}-h/2,h/2{)},\end{cases}
    }
$$ {of the interval $(-h/2,h/2)$} by trigonometric polynomials. In
the paper \cite{bk1} the sharp values of the best approximation
for  the special values of $h$ were found. In \cite{bk2} we gave
the complete solution of the problem for arbitrary values of
$h{\in(0,1]}$. In general case \cite{bk2} the situation is more
deep and results are not so simple as in \cite{bk1}.  For
applications to the  problem of optimal constants in the
Jackson-type inequalities we need, however,  results on
$L${-}approximation of $B${-}splines and linear combinations of
$B${-}splines (see \cite{fks,bk3}). Here we present some simple
results about $L${-}approximation of $B${-}splines as well as give
the the proof of its sharpness for the special values of $h$. In
some sense we give the appendix to the paper \cite{bk1}. \vskip
.4cm

The $B${-}splines are the convolutions of function  $\chi_h$ with
{itself}: $$ \chi_h^1(x):=\chi_h(x), $$

$$
  \chi_{h}^k (x):=\int_{{\mathbb{R}}}\chi_h(u)\chi_h^{k-1}(x-u)\,du=(\chi_h*\chi_h^{k-1})(x).
$$
\noindent
The $B$-splines are the functions with the $k-1$ order smoothness  and the supports
$
\mbox{supp}\ \chi_h^k={(}-kh/2, kh/2{)},\quad|\mbox{supp}\
\chi_{h}^k|=kh.
$

\noindent
Particularly

$$ \chi_{h}^2 (x)= \begin{cases} h^{-1} (1-|x|h^{-1}), \quad &
x\in{(}-h,h{)}, \\
 0, \quad & x \notin{(}-h,h{)}. \end{cases}
$$ \noindent It is easy to check that the operator of $k${-}th
order differentiation transforms $k${-}th $B${-}splines to
$k${-}th central differences:

$$
D^k (f*\chi_{h}^k) (x)= h^{-k}  \Delta_{h}^k f(x),
$$

$$
{\Delta}_{h}^k f(x) := \sum_{j=0}^k (-1)^j \binom{k}j f(x+kh/2
-jh).
$$
\noindent
Since

$$
\sup_{x} | \Delta_h^k f (x)| = \|{\Delta}_{h}^k f \| \le \sum_{j=0}^k \binom k j \| f \| \le 2^k \| f
\|,
$$
then
$$
\| D^k (f*\chi_{h}^k) \| \le (h/2)^{-k} \| f \|.  \eqno (1)
$$

\noindent
One of the main tools in  approximation theory is
the classical Favard's \cite{f} inequality:

$$ \| g \| \le \mathcal{K}_k (2\pi n)^{-k} \| D^k g  \|, \quad g \in T_{2n-1}^\perp,
\quad \mathcal{K}_k:= \frac 4 \pi \sum_{j= -\infty}^{\infty} \frac
1{(4j+1)^{k+1}}. \eqno (2)$$
\noindent
The Favard's constants $\mathcal{K}_k$  have the following properties

$$
1=\mathcal{K}_0 < \mathcal{K}_2=\pi^2/8 < \cdots < 4/\pi < \cdots < \mathcal{K}_3=\pi^3/24<
\mathcal{K}_1=\pi/2.
$$
\noindent
Direct consequence of (2) and (1) is

$$
\| g*\chi_{h}^k \| \le \mathcal{K}_k (2 \pi n)^{-k} \|D^k ( g * \chi_{h}^k) \| \le
\mathcal{K}_k (\pi nh)^{-k} \| g \|, \quad g \in T_{2n-1}^\perp.
$$
\vskip .3cm
\noindent
Therefore, we have
$$
E_{n} (\chi_{h}^k)_1 =
\sup_{g \in T_{2n-1}^\perp , \ \|g \| =1}
\int_T g(u) \chi_{h}^k (- u) \, du
$$
$$
\le \sup_{g \in T_{2n-1}^\perp , \ \|g \|
=1}| (\chi_{h}^k *g )(0)|
\le \mathcal{K}_k (\pi nh)^{-k}. \eqno (3)
$$
\break

{\bf Theorem. } {\it Let\quad {$k,n\in\mathbb{N}$,\quad
$k\le{n},$}\quad $h(\alpha)=\alpha/(2n),\quad 0 < \alpha \,
{\le}\, 2n/k $.\quad Then

$$ E_{n} (\chi_{h(\alpha)}^k )_1  \le  F_k\alpha^{-k},  \quad
\mbox{where} \  \quad F_k:=(2/\pi)^k \mathcal{K}_k. \eqno (4) $$
\vskip .2cm

\noindent
For example for  $k=1,2,3 $
we have

$$
E_{n} (\chi_{h(\alpha)} )_1  \le \frac 1{\alpha},
$$

$$
E_{n} (\chi^2_{h(\alpha)})_1  \le \frac{ 1} {2\alpha^2},
$$

$$
E_{n} (\chi^3_{h(\alpha)} )_1  \le \frac{1} {3  \alpha^3}.
$$

The inequalities {$(4)$} become equalities if $\alpha = 2j+1, \
{j\in\mathbb{Z}_+, \ j\le\displaystyle\frac{2n-k}{2k}}$. }

\vskip .2cm

The question about the value of the best $L${-}approximation of
$B${-}spline for {\it arbitrary} $0<\alpha \,{\leq} \, 2n/k$ is
not so simple (see \cite{bk2} for the case $k=1$) .

\vskip .2cm
{\it Proof.} The estimate (4) follows from the inequality (3).
We need to prove equalities for $\alpha = 2j+1$ only. At first consider the case
$k=1$.  We will use notation
 $$
 c_y(x):=\cos (2\pi x y), \quad y\in\mathbb{R}.
 $$
\noindent The function $ \pm \, \mbox{sign} \,(c_n),$
${n\in\mathbb{N}}$ gives  equality in (4). In other words, for
$k=1$ and \mbox{$h_j=(2j+1)/(2n),\  {j=0, \dots, n-1}$} we have $$
E_{n}(\chi_{h_j})_1\ge
\int_{{\mathbb{R}}}\chi_{h_j}(u)(-1)^j\mbox{sign}\,c_n(u)\,du=1/(2j+1).\eqno(5)
$$ \noindent One can rewrite the equality in (5) as

$$
\int_{{\mathbb{R}}}\chi_{2j+1}(u)(-1)^j\mbox{sign}\,c_{1/2}(u)\,du=1/(2j+1).\eqno
(6) $$
\noindent
 Note, that  $\mbox{sign} \,( c_{1/2} (x))\equiv
\mathcal{E}_0(x)$, where $\mathcal{E}_0(x)$ is the first Euler's
spline (see \cite{dl}, pp. 148--151). The Euler splines
$\mathcal{E}_k(x)$ are defined as follows:

$$
 \mathcal{E}_{j+1}(t)=\gamma_j\int_{{\mathbb{T}}}\mathcal{E}_{j}(x+u)\,du,
 \quad\gamma_j^{-1}=\int_{-1/2}^{1/2}\mathcal{E}_{j}(u)\,du,
$$ \noindent and have the following properties:
 $$
 \mathcal{E}_{j}(x+2)= \mathcal{E}_{j}(x), \quad  \mathcal{E}_{j}(x+1)=
 -\mathcal{E}_{j}(x),
$$

$$
\int_{-1}^{1} \mathcal{E}_{j}(u+x) \, du =0,
$$

$$ \mathcal{E}_{j}(-x) = \mathcal{E}_{j}(x), \quad
\mathcal{E}_{j}(-x-1/2) = \mathcal{E}_{j}(x+1/2), $$

$$ \| \mathcal{E}_{j} \| =1, \quad  \mathcal{E}_{j}(\nu) =
(-1)^\nu, \quad\nu\in{\mathbb{N}}, $$

$$
 D \mathcal{E}_{j}(x) = \pi \mathcal{K}_{j-1} \mathcal{K}_j^{-1}
\mathcal{E}_{j-1}(x+1/2). \eqno (7) $$

Come back to  $(6)$. Integrating by parts (7) we get

$$ (-1)^j\int_{{\mathbb{R}}}\chi_{2j+1} (u) \mathcal{E}_{0}(u)\,du
= $$ $$ \frac{(-1)^{j}}{2}\int_{{\mathbb{R}}}\chi_{2j+1} (u)\
D\mathcal{E}_{1}(u-1/2)\,du
=
\frac{(-1)^{j+1}}{2}\int_{{\mathbb{R}}}D \chi_{2j+1} (u) \
\mathcal{E}_{1}(u-1/2)\,du= $$

$$
\frac{(-1)^{j+1}}{2}(2j+1)^{-1}\int_{{\mathbb{R}}}\Delta_{2j+1}^1\delta(u)\
\mathcal{E}_{1}(u-1/2) \,du=\frac{(-1)^{j+1}}{2} (2j+1)^{-1}
(-\Delta_{2j+1}^1 \ \mathcal{E}_{1}(-1/2))= $$

$$
\frac{(-1)^{j+1}}{2}(2j+1)^{-1}(\mathcal{E}_{1}(-j-1)-\mathcal{E}_{1}(j))=
$$

$$ \frac{(-1)^{j+1}}{2}(2j+1)^{-1}((-1)^{j+1}-(-1)^j) =
(2j+1)^{-1}. $$
 \noindent One can rewrite the proof for odd $k$
without essential modifications: $$
\int_{{\mathbb{R}}}\chi_{h_j}^k(u)(-1)^j\mbox{sign}\,c_n(u)\,du=
\int_{{\mathbb{R}}}\chi_{2j+1}^k(u)(-1)^j\mbox{sign}\,c_{1/2}(u)\,du=
$$

$$
\frac{(-1)^j\mathcal{K}_k}{\pi^k}\int_{{\mathbb{R}}}\chi_{2j+1}^k(u)\
D^k\mathcal{E}_k
(u-k/2)\,du=\frac{(-1)^{j+1}\mathcal{K}_k}{\pi^k}\int_{{\mathbb{R}}}D^k\chi_{2j+1}^k(u)\
\mathcal{E}_k (u-k/2)\,du= $$

$$ \frac{(-1)^{j+1}\mathcal{K}_k}{\pi^k}
(2j+1)^{-k}\int_{{\mathbb{R}}}\Delta_{2j+1}^k\delta(u)\
\mathcal{E}_k (u-k/2)\,du= $$

$$
\frac{(-1)^{j+1}\mathcal{K}_k}{\pi^k}(2j+1)^{-k}(-\widehat\Delta_{2j+1}^k\mathcal{E}_k(-k/2))=
\frac{(-1)^{j+1}\mathcal{K}_k}{\pi^k}(2j+1)^{-k} 2^k(-1)^{j+1}=\frac{F_k}{(2j+1)^{k}}.
$$
\noindent
 Consider the case of even $k$.
\noindent We will use the equality $$ D^k \mathcal{E}_k (x) =
\frac {\pi^k}{\mathcal{K}_k} (-1)^{k/2} \mathcal{E}_0 (x), $$
\noindent which  implies $$ \int_{{\mathbb{R}}}\chi_{2j+1}^k(u)\
\mbox{sign}\,(c_{1/2}(u))\,dt=
(-1)^{k/2}\frac{\mathcal{K}_k}{\pi^k}\int_{{\mathbb{R}}}\chi_{2j+1}^k(u)\
D^k\mathcal{E}_k (u)\,du. $$ \noindent The integration by parts
gives $$ \int_{{\mathbb{R}}}\chi_{2j+1}^k(u)\ D^k \mathcal{E}_k
(u)\,du= \int_{{\mathbb{R}}}D^k\chi_{2j+1}^k(u)\ \mathcal{E}_k
(u)\,du. $$ Since $$
D^k\chi_{2j+1}^k(u)=(2j+1)^{-k}\Delta_{2j+1}^k\delta(u), $$ then
$$
(2j+1)^k\int_{{\mathbb{R}}}D^k\chi_{2j+1}^k(u)\mathcal{E}_k(u)\,du=
(\Delta_{2j+1}^k\delta*\mathcal{E}_k
)(0)=\Delta_{2j+1}^k\mathcal{E}_k(0)=2^k (-1)^{k/2}.\qed $$ \vskip
0.2cm {\bf Remark  1.} The restriction  $\alpha k \le 2n$ in the
Theorem means that we work with usual $B$-splines with support in
$[-1/2, 1/2]$. The inequality (4) is true without this
restriction, but we do not consider the sharpness of (4) for other
values of $\alpha k$. \vskip .2cm {\bf Remark 2.} The
approximations of  characteristic functions $(k=1)$ were
considered with details in \cite{bk1,bk2}. Note, that  one can
obtain the nontrivial approximation of step--function iff  $\alpha
> 1$. For small $\alpha$ the polynomial of the best approximation
of the step--function is equal to $0$, and the best approximation
is equal to $1$. As in the case of step--function we can indicate
the value of parameter $\alpha $ for the nontrivial estimates of
the best approximation. But we do not know what is the critical
value of $\alpha_0$ for nontrivial approximation if $k>1$. For $
\alpha > k^{-1/k}, \  k=2,3, \ $ we have $ E_n
(\chi_{h(\alpha)}^k)_1 < 1$ but we do not know the best
approximation in the case $\alpha = k^{-1/k}$. The $x_{minimum}$
of the function  $x^{-1/x}$ lies in $[2,3]$. Probably, there are
some links between this fact and the optimal smoothness of the
averaging operators. The averaging of the second order (
convolution with the hat function $\chi_h^2$ ) often gives the
most useful and sharp results.

\vskip .2cm

{\bf Remark 3.} This remark is close to { Remark 2}. We gave the
Theorem in { simple form}. We can present  here the more precise
version of (4):

$$
E_{n} (\chi_{h(\alpha)}^k )_1  \le  \min \left( 1, \frac{ F_k } { \alpha^k} \right).
\eqno (4')
$$

Note, that  for  $ \alpha  \le  F_k^{1/k}$ the inequality ($4'$)
gives trivial estimate. For  $ \alpha  > F_k^{1/k}$  (in other
words, if the support of the $k${-}th
 $B${-}spline  $\chi_{h}^k$ is greater then $ k
F_k^{1/k}/(2 n)$  )  the best approximation of $\chi_h^k$ is less then $1$.

\vskip .2cm
\small
\begin{tabular}[t]{l}
A.G. Babenko \\ Institute of Mathematics and Mechanics \\ Ural
Branch of the Russian Academy of Sciences \\ 16, S.Kovalevskoi
Str. \\ Ekaterinburg, 620219 \\ Russia \\
\textrm{babenko@imm.uran.ru}
\end{tabular}
 \hfill
\vskip .3cm
\begin{tabular}[t]{l}
Yuri Kryakin \\ Institute of Mathematics \\ University of
Wroc{l}aw
\\ Plac Grunwaldzki 2/4 \\ 50-384 Wroc{l}aw \\ Poland \\
\textrm{kryakin@gmail.com}
\end{tabular}
\end{document}